\documentclass{amsart}
\RequirePackage[pdftex=true,twoside=true,paperwidth=156mm,paperheight=234mm,bottom=25mm,left=23mm,right=18mm,headsep=7mm,headheight=8pt,top=29mm]{geometry} 

\usepackage{hyperref}

\usepackage{amsmath}
\usepackage{amssymb}
\usepackage{amsthm}

\usepackage{tikz}
\usetikzlibrary{fadings}

\def\AA{{\mathcal A}}

\def\CC{{\mathcal C}}

\def\IN{{\Bbb{N}}}

\def\IQ{{\Bbb{Q}}}
\def\IR{{\Bbb{R}}}

\def\TO{\Longrightarrow}
\def\In{\subseteq}

\def\prefix{\sqsubseteq}

\def\mto{\rightrightarrows}

\def\J{{\mathsf{{J}}}}

\def\L{{\mathsf{{L}}}}

\def\C{{\mathsf{{C}}}}

\def\Tr{{\mathrm{Tr}}}

\def\ACC{\mathsf{ACC}}
\def\DNC{\mathsf{DNC}}
\def\ATR{\mathsf{ATR}}
\def\LPO{\mathsf{LPO}}
\def\LLPO{\mathsf{LLPO}}
\def\MLPO{\mathsf{MLPO}}
\def\WKL{\mathsf{WKL}}
\def\RCA{\mathsf{RCA}}
\def\ACA{\mathsf{ACA}}
\def\SEP{\mathsf{SEP}}

\def\C{\mathsf{C}}

\def\SORT{\mathsf{SORT}}
\def\RAT{\mathsf{RAT}}
\def\SEP{\mathsf{SEP}}

\def\MCT{\mathsf{MCT}}

\def\leqW{\mathop{\leq_{\mathrm{W}}}}
\def\equivW{\mathop{\equiv_{\mathrm{W}}}}

\def\leqSW{\mathop{\leq_{\mathrm{sW}}}}
\def\equivSW{\mathop{\equiv_{\mathrm{sW}}}}

\def\nleqW{\mathop{\not\leq_{\mathrm{W}}}}

\def\lW{\mathop{<_{\mathrm{W}}}}
\def\lSW{\mathop{<_{\mathrm{sW}}}}

\def\bigtimes{\mathop{\mathsf{X}}}

\def\leqW{\mathop{\leq_{\mathrm{W}}}}
\def\equivW{\mathop{\equiv_{\mathrm{W}}}}

\def\leqSW{\mathop{\leq_{\mathrm{sW}}}}
\def\equivSW{\mathop{\equiv_{\mathrm{sW}}}}

\def\nleqW{\mathop{\not\leq_{\mathrm{W}}}}

\def\lW{\mathop{<_{\mathrm{W}}}}
\def\lSW{\mathop{<_{\mathrm{sW}}}}

\def\id{{\mathrm{id}}}

\def\dom{{\mathrm{dom}}}
\def\range{{\mathrm{range}}}

\def\NON{\mathsf{NON}}

\def\DIS{\mathsf{DIS}}

\def\EC{\mathsf{EC}}

\newcommand{\uwidehat}[1]{%
  \mathpalette\douwidehat{#1}%
}
\makeatletter
\newcommand{\douwidehat}[2]{%
  \sbox0{$\m@th#1\widehat{\hphantom{#2}}$}%
  \sbox2{$\m@th#1x$}
  \sbox4{$\m@th#1#2$}
  \dimen0=\ht0
  \advance\dimen0 -.8\ht2
  \dimen2=\dp4
  \rlap{%
    \raisebox{\dimexpr\dimen0-\dimen2}{%
      \scalebox{1}[-1]{\box0}%
    }%
  }%
  {#2}%
}
\makeatother

\usepackage{catchfile}

\title[Weihrauch Complexity and the Hagen School]{Weihrauch Complexity and the\\ Hagen School of Computable Analysis}

\author{Vasco Brattka}
\address{Faculty of Computer Science, Universit\"at der Bundeswehr M\"unchen, Germany and
Department of Mathematics and Applied Mathematics, University of Cape Town, South Africa}
\email{Vasco.Brattka@cca-net.de}

\begin{document}

\theoremstyle{definition}
\newtheorem{theorem}{Theorem}
\newtheorem{definition}[theorem]{Definition}
\newtheorem{problem}[theorem]{Problem}
\newtheorem{assumption}[theorem]{Assumption}
\newtheorem{corollary}[theorem]{Corollary}
\newtheorem{proposition}[theorem]{Proposition}
\newtheorem{lemma}[theorem]{Lemma}
\newtheorem{observation}[theorem]{Observation}
\newtheorem{fact}[theorem]{Fact}
\newtheorem{question}[theorem]{Question}
\newtheorem{example}[theorem]{Example}
\newtheorem{convention}[theorem]{Convention}
\newtheorem{conjecture}[theorem]{Conjecture}


\begin{abstract}
Weihrauch complexity is now an established and active part of mathematical logic.
It can be seen as a computability-theoretic approach
to classifying the uniform computational content of mathematical problems.
This theory has become an important interface between more proof-theoretic and more computability-theoretic 
studies in the realm of reverse mathematics.
Here we present a historical account of the early developments of 
Weihrauch complexity by the Hagen school of computable analysis
that started more than thirty years ago, and we indicate how this has influenced, informed,
and anticipated more recent developments of the subject. 
\end{abstract}

\maketitle

\setcounter{tocdepth}{1}
\tableofcontents

\section{Computable Analysis in Hagen}

The Hagen school of computable analysis was founded by Klaus Weihrauch in the 1980s. 
Klaus Weihrauch received his PhD from the University of Bonn in 1973, and after a research associateship
at Cornell University and a professorship position at RWTH Aachen he held a
chair for Theoretical Computer Science at the University of Hagen from 1979 until his retirement in 2008.
Together with his PhD students he developed the representation
based approach to computable analysis, sometimes called {\em type-2 theory of effectivity}. 
Most notably among the early PhD students were Christoph Kreitz~\cite{Kre84}, Thomas Deil~\cite{Dei84a}
and Norbert M\"uller~\cite{Mue88}. A later generation of PhD students includes Peter Hertling~\cite{Her96d}, 
Xizhong Zheng~\cite{Zhe98a}, Matthias Schr\"oder~\cite{Sch02}, 
 the author of these notes and others.\footnote{A more complete PhD genealogy of Klaus Weihrauch can be found in 
 the preface of~\cite{BHKZ04}.}

The main idea of this approach to computable analysis is to perform all computability  considerations on 
Baire space, and to transfer theses concepts to other spaces by representing them with Baire space.
On the one hand, on Baire space $\IN^\IN$ concepts such as continuity and computability are 
well understood, for instance with the help of Turing machines that operate on natural number sequences. 
On the other hand, natural number sequences can naturally be used to represent other objects
such as real numbers, closed subsets or continuous functions on real numbers.
Analogous statements hold for Cantor space $2^\IN$, which additionally has natural concepts of time and space complexity.

Implicitly, such representations were used ever since Turing~\cite{Tur37,Tur38a} introduced his machines to operate on real numbers,
and they are also implicit in constructive analysis~\cite{Bis67}, reverse mathematics~\cite{Sim99},
descriptive set theory~\cite{Mos80,Kec95} and in set-theoretical constructions in classical mathematics too. 
However, the crucial idea of Klaus Weihrauch and his collaborators was to make such representations
objects of mathematical investigations themselves by considering them as partial surjective maps.

\begin{definition}[Representation]
\label{def:representation}
A {\em representation} of a set $X$ is a partial surjective map $\delta:\In\IN^\IN\to X$.
\end{definition}

This idea already had a tradition in mathematical logic since it was also present
in Hauck's work~\cite{Hau78,Hau80a}.  
Kreitz and Weihrauch~\cite{KW85,Wei85,KW87,WK87} developed the theory of representations following the lines
of the theory of numberings that was proposed in the 1970s by Ershov~\cite{Ers99}.
An early manifesto by Kreitz and Weihrauch can be found in~\cite{KW84}
and more complete presentations in~\cite{Wei87,Wei00}.

It turned out that when dealing with infinite objects such as real numbers
the proper choice of a representation is more crucial than for discrete objects.
When one represents rational numbers $\IQ$ one more or less automatically arrives
at a suitable representation from the perspective of computability,
and one has to work hard to find a representation that is not computably equivalent to a natural one, 
for instance by artificially encoding the halting problem into the representation.

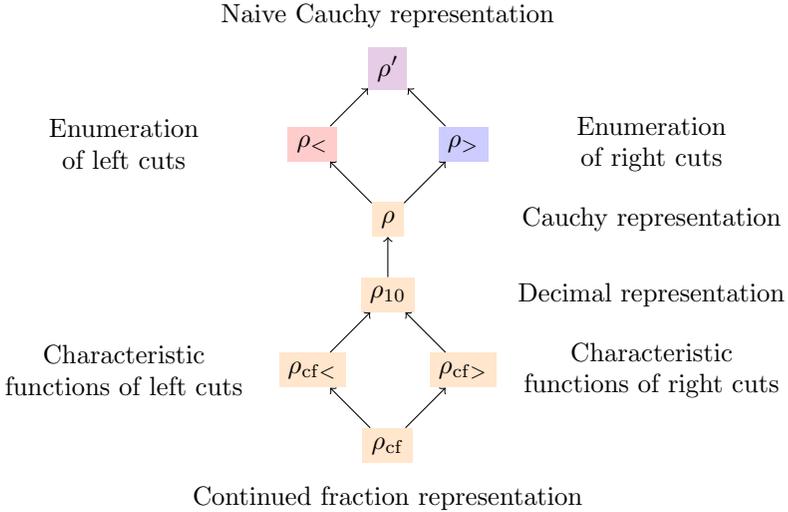
\begin{figure}[htb]
\begin{tikzpicture}[scale=1,auto=left]
\node[style={fill=orange!20}]  (RCF) at (0,0) {$\rho_{\mathrm{cf}}$};
\node[style={fill=orange!20}]  (RCF<) at (-1,1) {$\rho_{\mathrm{cf}<}$};
\node[style={fill=orange!20}]  (RCF>) at (1,1) {$\rho_{\mathrm{cf}>}$};
\node[style={fill=orange!20}]  (R10) at (0,2) {$\rho_{10}$};
\node[style={fill=orange!20}]  (R) at (0,3) {$\rho$};
\node[style={fill=red!20}]  (R<) at (-1,4) {$\rho_{<}$};
\node[style={fill=blue!20}]  (R>) at (1,4) {$\rho_{>}$};
\node[style={fill=violet!20}]  (RP) at (0,5) {$\rho'$};

\draw[->] (RCF) edge (RCF<);
\draw[->] (RCF) edge (RCF>);
\draw[->] (RCF<) edge (R10);
\draw[->] (RCF>) edge (R10);
\draw[->] (R10) edge (R);
\draw[->] (R) edge (R<);
\draw[->] (R) edge (R>);
\draw[->] (R<) edge (RP);
\draw[->] (R>) edge (RP);

\node at (0,5.7) {Naive Cauchy representation};
\node at (3.5,4) {\begin{minipage}{3cm}\begin{center}Enumeration\\ of right cuts\end{center}\end{minipage}};
\node at (-3.5,4) {\begin{minipage}{3cm}\begin{center}Enumeration\\ of left cuts\end{center}\end{minipage}};
\node at (3.5,3) {Cauchy representation};
\node at (3.5,2) {Decimal representation};
\node at (3.5,1) {\begin{minipage}{4cm}\begin{center}Characteristic\\functions of right cuts\end{center}\end{minipage}};
\node at (-3.5,1) {\begin{minipage}{4cm}\begin{center}Characteristic\\functions of left cuts\end{center}\end{minipage}};
\node at (0,-0.7) {Continued fraction representation};
\end{tikzpicture}
\caption{Lattice of real number representations.}
\label{fig:real-number-representations}
\end{figure}

In the case of infinite objects, such as real numbers $\IR$, there are many natural representations
that lead to mutually inequivalent structures. Diagram~\ref{fig:real-number-representations} displays a portion
of the lattice of real number representations that were studied already by Deil~\cite{Dei84a}.
We omit the formal definitions, but we introduce
some symbolic names in the diagram for later reference.
The representations of the same color induce identical notions of computable real numbers: the Cauchy representation
and all representations below it in the diagram induce the ordinary notion of a computable real number; the representations
via enumerations of left and right cuts induce the left- and right-computable real numbers, respectively; the native
Cauchy representation induces the limit computable reals. 

The preorder used in the diagram is that of reducibility of representations,
which can be seen as a generalization of the concept of reducibility for numberings
or of many-one reducibility to type-2 objects.
Any arrow in the diagram indicates a reduction in the direction of the arrow,
and all missing arrows indicate that no reduction is possible (except for those that follow by transitivity and reflexivity).

\begin{definition}[Reducibility]
\label{def:reducibility}
Let $f,g:\In\IN^\IN\to X$ be partial functions. Then we say that $f$ {\em is reducible to} $g$,
in symbols $f\leq g$, if there exists a computable $F:\In\IN^\IN\to\IN^\IN$ such that
$f(p)=g\circ F(p)$ for all $p\in\dom(f)$.
\end{definition}

By $\equiv$ we denote the equivalence induced by $\leq$. 
One can also define a purely topological analogue of this reducibility by requiring that
$F$ is continuous, and it turns out that the relations in the diagram in Figure~\ref{fig:real-number-representations}
 are not affected by this modification. Hence, one can say that the uniform distinctions between these
 representations are already of topological nature.
 
In particular, the continued fraction representation carries more continuously accessible information
 about real numbers than the decimal representation, and in turn the decimal representation
 carries more information than the Cauchy representation.
 Having more informative representations on the input side is helpful, but it can be a burden
 on the output side. For instance addition is neither computable with respect to the continued fraction
 representation~\cite{Ko86b} nor with respect to the decimal representation~\cite{Tur38a}.

These observations naturally lead to the question how one can identify a suitable representation
among all the many representations of infinite objects that one can consider?
The answer that Kreitz and Weihrauch gave is that a good representation has to be
topologically natural, and the corresponding concept is called {\em admissibility}.
Their main theorem on this topic is the following~\cite{KW84,KW85,Wei87}.

\begin{theorem}[Kreitz-Weihrauch 1984]
\label{thm:Kreitz-Weihrauch}
If $X$ and $Y$ are admissibly represented $T_0$--spaces with countable bases,
then $f:\In X\to Y$ is continuous if and only if it is continuous with respect to the underlying representations.
\end{theorem}

Some further explanations are required here. For one, continuity with respect to the underlying representations
means that the diagram in Figure~\ref{fig:admissibility} commutes.

\begin{figure}[htb]
\begin{center}
\unitlength 0.60mm
\linethickness{0.4pt}
\begin{picture}(59.00,59.00)
(0,5)
\put(5.00,55.00){\makebox(0,0)[cc]{$\IN^\IN$}}
\put(10.00,55.00){\vector(1,0){40.00}}
\put(5.00,5.00){\makebox(0,0)[cc]{$X$}}
\put(10.00,5.00){\vector(1,0){40.00}}
\put(30.00,59.00){\makebox(0,0)[cc]{$F$}}
\put(30.00,1.00){\makebox(0,0)[cc]{$f$}}
\put(5.00,50.00){\vector(0,-1){40.00}}
\put(59.00,30.00){\makebox(0,0)[lc]{$\delta_Y$}}
\put(1.00,30.00){\makebox(0,0)[rc]{$\delta_X$}}
\put(55.00,55.00){\makebox(0,0)[cc]{$\IN^\IN$}}
\put(55.00,5.00){\makebox(0,0)[cc]{$Y$}}
\put(55.00,50.00){\vector(0,-1){40.00}}
\end{picture}
\end{center}
\caption{Continuity with respect to representations.}
\label{fig:admissibility}
\end{figure}

That is, if $(X,\delta_X)$ and $(Y,\delta_Y)$ are represented spaces, then a function $f:\In X\to Y$ is called 
{\em continuous with respect to the underlying representations}, if there is is a continuous
$F:\In\IN^\IN\to\IN^\IN$ such that $\delta_YF(p)=f\delta_X(p)$ for all $p\in\dom(f\delta_X)$.
Other notions such as computability, Borel measurability, etc., can be transferred analogously
to represented spaces, and Theorem~\ref{thm:Kreitz-Weihrauch} tells us that as long we
use admissible representations, then at least topological continuity is the same as continuity with respect 
to the representations. 

Theorem~\ref{thm:Kreitz-Weihrauch} was later generalized by Schr\"oder~\cite{Sch02} to a larger category of
topological spaces, and one needs to replace topological continuity with sequential continuity
in this more general context. Schr\"oder also introduced a more general definition of admissibility
that we are going to present here.

\begin{definition}[Admissibility]
A representation $\delta:\In \IN^\IN\to X$ of a topological space $X$ is called 
{\em admissible} if it is continuous
and maximal among all continuous representations of $X$ with respect to the topological version
of the reducibility $\leq$.
\end{definition}

\begin{figure}[htb]
\includegraphics[width=11cm]{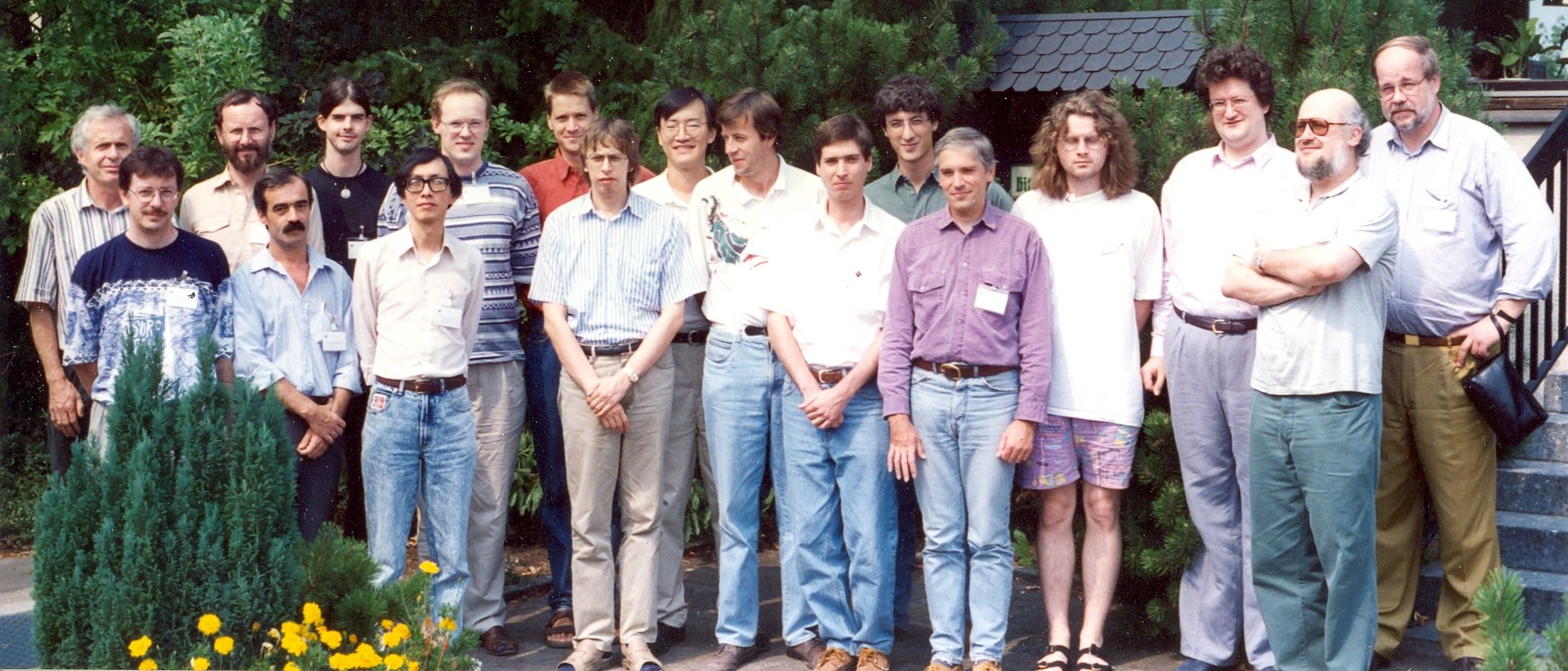}
\caption{A group picture with some participants of CCA 1995 (taken by the author). Names from left to right:
Klaus Weihrauch, Norbert M\"uller, Ludwig Staiger, Kostas Skandalis, Markus Bl\"aser,
Ker-I Ko, Peter Hertling, Jens Blanck, Matthias Schr\"oder, Arthur Chou, Klaus Meer, 
Martin H\"otzel Escardo, Pietro Di Gianantonio, Uwe Mylatz, Holger Schulz, Rudolf Freund, Janos Blazi, Achim Kallweit.}
\label{fig:CCA1995}
\end{figure}

Among the real number representations it is the equivalence class of the Cauchy representation
that yields admissible representations with respect to the Euclidean topology.
In general, admissibility yields a handy criterion to judge whether
a representation is suitable from a topological perspective.
Many hyper and function space representations have been analyzed, such as representations for the space $\CC(X)$ of continuous functions
$f:X\to \IR$, which we are not going to introduce in detail here~\cite{BHW08,Wei00}.
There are also other schools of computable analysis that approach the subject from a slightly different angle,
such as
the Pour-El and Richards school~\cite{PR89}. 
Studies of computational complexity in analysis were initiated by Norbert M\"uller~\cite{Mue86a,Mue93,MM93}, Ker-I Ko~\cite{Ko91}, and more
recently by Akitoshi Kawamura and Stephen Cook~\cite{Kaw10,KC12}. 
A more comprehensive discussion of the
historical developments in computable analysis that arose out of Turing's work is presented in~\cite{AB14}.
The tutorial~\cite{BHW08} contains a concise introduction to computable analysis, and
the handbook~\cite{BH21} contains surveys on many fascinating aspects of the more recent developments in computable analysis in general.
From 1995 until today the computable analysis group founded in Hagen runs 
a conference with the title {\em Computability and Complexity in Analysis} (CCA).
Some participants of the first meeting CCA 1995 in Hagen are shown in Figure~\ref{fig:CCA1995}.

\section{Weihrauch Reducibility}

When studies in computable analysis progressed, the need arose to include the study of multi-valued
maps. For instance, when one discusses the problem of finding zeros of a continuous
function $f:X\to\IR$, then it is not sufficient to consider single-valued functions of type $Z:\In\CC(X)\to\IR$
in oder to describe zero finding,
because under certain assumptions it might be possible to compute a zero of a continuous
function $f:X\to\IR$ only in a non-extensional way, i.e., such that the zero depends on the
description of $f$ and not just on $f$ itself~\cite{Wei00}. This phenomenon is best captured by 
describing zero finding as a multi-valued partial function of type $Z:\In\CC(X)\mto\IR$. In general, many mathematical problems
can be construed as such multi-valued maps. Hence we use the following general definition.

\begin{definition}[Problem]
A {\em problem} $f:\In X\mto Y$ is a partial multi-valued map on represented spaces $X,Y$.
\end{definition}

If we have two problems $f,g:\In X\mto Y$ of the same type, then we say that
$f$ {\em refines} $g$, in symbols $f\prefix g$, if $\dom(g)\In\dom(f)$ and  $f(x)\In g(x)$ for all $x\in\dom(g)$.
Problems can naturally be combined in several ways. For instance, the {\em composition} $g\circ f$ of 
two problems $f:\In X\mto Y$ and $g:\In Y\mto Z$ is defined by $\dom(g\circ f):=\{x\in \dom(f):f(x)\In\dom(g)\}$ and
\[(g\circ f)(x):=\{z\in Z:(\exists y\in f(x))\;z\in g(y)\}\]
for all $x\in\dom(g\circ f)$.
The particular choice of the domain is important,
as this ensures that composition preserves computability and continuity in the appropriate way.
Another operation on the problems $f:\In X\mto Y$ and $g:\In W\mto Z$ is the {\em product} $f\times g$
that is defined by $\dom(f\times g):=\dom(f)\times\dom(g)$ and
\[(f\times g)(x,w):=f(x)\times g(w)\]
for all $(x,w)\in\dom(f\times g)$.

As soon as problems are captured as multi-valued maps, there arises the need for a tool to compare
the computational power of problems beyond refinement.
In some sense, the concept of reducibility given in Definition~\ref{def:reducibility} already gives us
a way to compare (single-valued) problems of a certain type. 
However, only the input is subject to a pre-processing step here.
Klaus Weihrauch's ideas of extending the notion of many-one reducibility such that also
a post-processing of the output is considered,
were  laid out in two unpublished technical reports:

\begin{itemize}
\item[\cite{Wei92a}] Klaus Weihrauch, {\em The degrees of discontinuity of some translators between representations of the real numbers}, Technical Report TR-92-050, International Computer Science Institute, Berkeley, July 1992.
\item[\cite{Wei92c}] Klaus Weihrauch, {\em The TTE-interpretation of three hierarchies of omniscience principles},
Informatik Berichte 130, FernUniversit\"at Hagen, September 1992.
\end{itemize}

In his original definition Klaus Weihrauch did not present the definition of his reducibilities
in the way it is seen most often now.
In order to make the definition as clear and simple as possible,
we present a slightly more general definition in almost categorical terms that can be interpreted in different ways.
By $\id:\IN^\IN\to\IN^\IN$ we denote the identity on Baire space.

\begin{definition}[Weihrauch reducibility]
Let $C$ be a class of problems, and let $f,g$ be two problems. We introduce the following terminology:
\begin{enumerate}
\item $f$ is {\em strongly Weihrauch reducible} to $g$ ({\em with respect to $C$}), in symbols $f\leqSW g$,
if there are $H,K\in C$ such that $H\circ g\circ K\prefix f$.
\item $f$ is {\em Weihrauch reducible} to $g$ ({\em with respect to $C$}), in symbols $f\leqW g$,
if there are $H,K\in C$ such that $H\circ(\id\times g)\circ K\prefix f$.
\end{enumerate}
The set $C$ has to be fixed in order to make the symbolic notation meaningful.
\end{definition}  

In the usual definition we use for $C$ the class of all computable (or continuous) problems.\footnote{That
this actually yields the ordinary definition of Weihrauch reducibility follows from \cite[Lemma~2.5]{BP18}.} 
And this is what usually is called (the topological version of) Weihrauch reducibility.
If $C$ is sufficiently nice, i.e., closed under composition and under product with $\id$, then
the above definitions actually yield preorders, no matter what $C$ is.
And in fact, versions of Weihrauch reducibility have also been considered for the classes $C$
of (suitably defined) polynomial-time computable problems~\cite{KC12},
arithmetic problems~\cite{MV21}, or even hyperarithmetic piecewise computable problems~\cite{GKT19}. 
In a more categorical setting concepts similar to Weihrauch reducibility were also studied by Hirsch~\cite{Hir90},
and in complexity theory similar concepts are known as (polynomial-time) many-one reducibility for functions on
discrete spaces~\cite{KC12}.

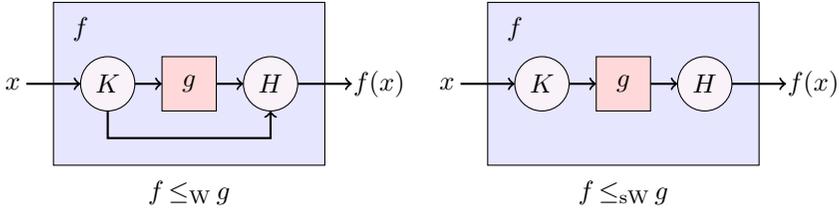
\begin{figure}[htb]
\begin{tikzpicture}[scale=.36,auto=left]
\draw[style={fill=blue!10}] (-3,6) rectangle (7,0);
\draw[style={fill=red!15}]  (1,4) rectangle (3,2);
\draw[style={fill=violet!5}]  (-1,3) ellipse (1 and 1);
\draw[style={fill=violet!5}]  (5,3) ellipse (1 and 1);
 \node at (-1,3) {$K$};
\node at (5,3) {$H$};
\node at (2,3) {$g$};
\node at (-2,5) {$f$};
\node at (-4.5,3) {$x$};
\node at (9,3) {$f(x)$};
\node at (2,-1) {$f\leqW g$};
\draw[->,thick] (0,3) -- (1,3);
\draw[->,thick] (-4,3) -- (-2,3);
\draw[->,thick] (3,3) -- (4,3) ;
\draw[->,thick] (6,3) -- (8,3);
\draw[->,thick] (-1,2) -- (-1,1) -- (5,1) -- (5,2);
\draw[style={fill=blue!10}] (13,6) rectangle (23,0);
\draw[style={fill=red!15}]  (17,4) rectangle (19,2);
\draw[style={fill=violet!5}]  (15,3) ellipse (1 and 1);
\draw[style={fill=violet!5}]  (21,3) ellipse (1 and 1);
 \node at (15,3) {$K$};
\node at (21,3) {$H$};
\node at (18,3) {$g$};
\node at (14,5) {$f$};
\node at (11.5,3) {$x$};
\node at (25,3) {$f(x)$};
\node at (18,-1) {$f\leqSW g$};
\draw[->,thick] (16,3) -- (17,3);
\draw[->,thick] (12,3) -- (14,3);
\draw[->,thick] (19,3) -- (20,3) ;
\draw[->,thick] (22,3) -- (24,3);
\end{tikzpicture}
\caption{Weihrauch reducibility and strong Weihrauch reducibility.} 
\label{fig:reducibility}
\end{figure}

The intuitive idea of (strong) Weihrauch reducibility is illustrated in the diagrams in Figure~\ref{fig:reducibility}.
The idea is that $f$ being {\em strongly Weihrauch reducible} to $g$ means that $g$ composed
with some suitable pre-processor $K$ and some suitable post-processor $H$ refines $f$.
Suitability includes that the pre- and post-processors have to be from the set $C$.
In the case of the ordinary (non-strong) reduction the post-processor additionally has access
to the original input in form of some information determined by the pre-processor. 

In the original definition in \cite{Wei92a,Wei92c} Weihrauch firstly only considered single-valued problems $f,g$ on Cantor space,
and for $C$ he used the set of continuous single-valued functions on Cantor space.
He denoted the reducibilities $\leqSW$ and $\leqW$ by $\leq_2$ and $\leq_1$, respectively (and he used $\leq_0$
for a many-one like reducibility as in Definition~\ref{def:reducibility}).
In a second step in \cite{Wei92a,Wei92c} he extended the definition to sets of such problems $f,g$.
But this is just another technical way of dealing with multi-valuedness and essentially
yields an approach that is equivalent to the modern one.\footnote{See the discussions in \cite[Section~2.1]{BP18} and \cite[Appendix~A]{DDH+16}.}
Hence, if not mentioned otherwise, we will assume from now on that $C$ is the class of all
computable problems, and we will express everything in modern terminology.
As usual, we denote the equivalences derived from strong and ordinary Weihrauch reducibility
by $\equivSW$ and $\equivW$, respectively, and the strict versions of the reducibilities
by $\lSW$ and $\lW$, respectively.

Only relatively late, it was discovered that the order structures induced by strong and ordinary Weihrauch reducibility
are lattices. For $\leqW$ this result is due to the author and Gherardi~\cite{BG11a} (for the lower semi-lattice) and Pauly~\cite{Pau10a} (for the upper semi-lattice).
For $\leqSW$ the result is due to Dzhafarov~\cite{Dzh19}.

\begin{proposition}[Weihrauch lattice]
The order structures induced by $\leqW$ and $\leqSW$ are lattices.
In the case of $\leqW$ the lattice is distributive, in the case of $\leqSW$ it is not.
\end{proposition}

More precise definitions along these lines and further results can be found in~\cite{BGP21}.

\section{Weihrauch Complexity}

One goal of Weihrauch complexity is to classify the computational content of mathematical theorems
of the logical form
\[(\forall x\in X)(x\in D\TO (\exists y\in Y)\;P(x,y))\]
in the Weihrauch lattice.
If $X$ and $Y$ are represented spaces, then such a theorem directly translates into a 
Skolem-like problem
\[F:\In X\mto Y,x\mapsto\{y\in Y:P(x,y)\}\]
with $\dom(F)=D$.
To locate the problem $F$ in the Weihrauch lattice amounts to classifying the uniform computational
complexity of the corresponding theorem. One way to calibrate the complexity of some $F$ is to compare
it to suitable {\em choice problems} for a suitable space $X$.
By $\C_X$ we denoted the {\em closed choice problem} 
\[\C_X:\In\AA_-(X)\mto X,A\mapsto A\]
for a computable metric space $X$.
The instances are non-empty closed sets $A\In X$ represented by negative information 
(the space of all such closed sets is denoted by $\AA_-(X)$)
and the solution can be any point in $A$ (see \cite{BBP12,BGP21} for more precise definitions).
Equivalently, we can see the instances as continuous functions $f:X\to\IR$ with zeros and the
solution can be any zero of the function $f$, i.e., if $\CC(X)$ denotes the space of continuous
function $f:X\to\IR$, represented in a natural way, then
\[\C_X:\In\CC(X)\mto X,f\mapsto f^{-1}\{0\},\]
with $\dom(\C_X)=\{f\in\CC(X):f^{-1}\{0\}\not=\emptyset\}$.
This alternative description of choice shows that choice problems are basically about
the solutions of equations of type
\[f(x)=0\]
for continuous $f:X\to\IR$.
Typical spaces $X$ whose choice problems have been considered are Cantor space $X=2^\IN$, Baire space $X=\IN^\IN$,
Euclidean space $X=\IR$, the natural numbers $X=\IN$, and finite spaces ${n=\{0,...,n-1\}}$ for $n\in\IN$.

\begin{figure}[htb]
\begin{tabular}{ll}
{\bf Weihrauch complexity} & {\bf reverse mathematics}\\\hline\\[-0.3cm]
$\C_1$ & $\RCA_0$ without $\Sigma^0_1$--induction\\
$\C_\IN$ & $\Sigma^0_1$--induction\\
$\C_{2^\IN}$ & $\WKL_0$ without $\Sigma^0_1$--induction\\
$\C_\IR$ & $\WKL_0$ with $\Sigma^0_1$--induction\\
$\C_{\IN^\IN}$ & $\ATR_0$
\end{tabular}
\caption{Choice problems versus axiom systems.}
\label{fig:reverse-mathematics}
\end{figure}

Proving equivalences to choice problems roughly corresponds to a uniform version of a  classification of the corresponding theorem
in reverse mathematics~\cite{Sim99}. Reverse mathematics is a proof-theoretic approach to classifying theorems according
to which axioms are needed to prove the theorem in second-order arithmetic. Typical axiom systems are {\em recursive comprehension} $\RCA_0$,
{\em arithmetic comprehension} $\ACA_0$, {\em arithmetic transitive recursion} $\ATR_0$. Details can be found in \cite{Sim99,Hir15}.
The table given in Figure~\ref{fig:reverse-mathematics} indicates the correspondences between Weihrauch complexity and reverse mathematics.

This is just a very rough correspondence. For instance, $\ATR_0$ can be analyzed more closely
in the vicinity of $\C_{\IN^\IN}$, and this is subject of several recent studies,
for instance, by Marcone and Valenti~\cite{MV21}, Goh, Pauly, and Valenti~\cite{GPV21}, Goh~\cite{Goh20a},
Kihara, Marcone and Pauly~\cite{KMP20}.

The system $\ACA_0$ can be characterized using iterations of the {\em limit problem} 
\[\lim:\In\IN^\IN\to\IN^\IN,\langle p_0,p_1,p_2,...\rangle\to\lim_{n\to\infty}p_n,\]
which is just the usual limit on Baire space (where for technical convenience, the input sequence is encoded by a standard tupling function $\langle\,\rangle$ in a single point in Baire space).
Using these choice problems, one can obtain the following results.
We are just mentioning some example for $\C_\IN$, $\C_{2^\IN}$ and $\lim$. Many further results can be found in the survey~\cite{BGP21}.
The following result on problems equivalent to $\C_\IN$ is due to the author and Gherardi~\cite{BG11a}.

\begin{theorem}[Choice on the natural numbers]
\label{thm:CN}
The following are all Weih\-rauch equivalent to each other:
\begin{enumerate}
\item Choice on natural numbers $\C_{\IN}$.
\item The Baire category theorem for computable complete metric spaces.
\item Banach's inverse mapping theorem for the Hilbert space $\ell_2$.
\item The open mapping theorem for $\ell_2$.
\item The closed graph theorem for $\ell_2$.
\item The uniform boundedness theorem on non-singleton computable Banach spaces.
\end{enumerate}
\end{theorem}

In order to be more precise, one would have to explain which logical version of the
theorem is translated into a problem here, and how the underlying data are represented.
That these aspects can make a difference has been discussed for the Baire category 
theorem in detail~\cite{BHK18}. We will not formalize these problems here and refer
the interested reader to~\cite{BG11a} for all details.

Next we mention a number of problems that are equivalent to $\C_{2^\IN}$.
The result on the Hahn-Banach theorem and the separation theorem is due to Gherardi and Marcone~\cite{GM09},
the result on the Brouwer fixed point theorem is due to the author, Le Roux, Joseph Miller,
and Pauly~\cite{BLRMP19}. The result on the Gale-Stewart theorem is due to Le Roux and Pauly~\cite{LRP15}.
All other results are easy to prove, the Heine-Borel theorem and the theorem of the maximum where briefly discussed in~\cite{Bra16}.

\begin{theorem}[Choice on Cantor space]
\label{thm:C2N}
The following are all Weihrauch equivalent to each other:
\begin{enumerate}
\item Choice on Cantor space $\C_{2^\IN}$.
\item Weak K\H{o}nig's lemma $\WKL$.
\item The Hahn-Banach theorem.
\item The separation theorem $\SEP$ on the separation of two disjoint enumerated sets in $\IN$ by the characteristic function of another set.
\item The Heine-Borel covering theorem.
\item The theorem of the maximum.\index{theorem of the maximum}
\item The Brouwer fixed point theorem for dimension $n\geq2$.
\item The theorem of Gale-Stewart (on determinacy of games on Cantor space with closed winning sets).
\end{enumerate}
\end{theorem}

Once again one would have to make these problems more precise. We give two examples. 
The {\em separation problem} is defined by
\[\SEP:\In\IN^\IN\mto2^\IN,\langle p,q\rangle\mapsto\{A\in 2^\IN:\range(p-1)\In A\In\IN\setminus\range(q-1)\}\]
with $\dom(\SEP):=\{\langle p,q\rangle:\range(p-1)\cap\range(q-1)=\emptyset\}$.
Here $p-1\in\IN^*\cup\IN^\IN$ is the finite or infinite sequence of natural numbers 
that consists of the concatenation of
\[p(0)-1,p(1)-1,p(2)-1,...\]
with the understanding that $-1$ is the empty word. This technical construction is used to allow for enumerations of
the empty set, and in order to keep $0$ as dummy value in enumerations.
We are going to see it later again.

Weak K\H{o}nig's lemma is defined by
\[\WKL:\In{\Tr}\mto2^\IN,T\mapsto[T],\]
where $\Tr$ is the set of binary trees, $\dom(\WKL)$ contains all infinite such trees,
and $[T]$ is the set of infinite paths of such a tree  $T$.
As a simple example we mention at least one proof~\cite[Proposition~2.8 and Theorem~2.11]{BG11a}.
If we represent the space $\AA_-(2^\IN)$ of closed subsets of $2^\IN$ by negative information,
then the map $f:\Tr\to\AA_-(2^\IN),T\to[T]$ is easily seen to be computable and it admits
a multi-valued computable right inverse $g:\AA_-(2^\IN)\mto\Tr$. This is all that is needed
to prove $\WKL\equivW\C_{2^\IN}$.

We close this section with a number of problems that are equivalent to the limit operation.
The results on the monotone convergence theorem and on the operator of differentiation $d$ (which is restricted to continuously differentiable functions)
are due to von Stein~\cite{Ste89} (see also Theorem~\ref{thm:von-Stein}).
The Radon-Nikodym theorem was studied by Hoyrup, Rojas, and Weihrauch~\cite{HRW12}.
The other results are easy to show (see for instance~\cite{BBP12,Bra16,Bra18}).

\begin{theorem}[The limit]
\label{thm:lim-thm}
The following are all Weihrauch equivalent to each other:
\begin{enumerate}
\item The limit map $\lim$ on Baire space (or Cantor space, or Euclidean space).
\item The Turing jump $\J:\IN^\IN\to\IN^\IN,p\mapsto p'$.
\item The monotone convergence theorem $\MCT:\In\IR^\IN\to\IR,(x_n)_{n\in\IN}\mapsto\sup_{n\in\IN}x_n$.
\item The operator of differentiation $d:\In\CC[0,1]\to\CC[0,1],f\mapsto f'$.
\item The Fr{\'e}chet-Riesz representation theorem for $\ell_2$.
\item The Radon-Nikodym theorem.
\end{enumerate}
\end{theorem}

In between the degrees of $\C_{2^\IN}$ and $\lim$ there is the degree of $\C_\IR$ and the degree
of the {\em lowness problem} $\L:=\J^{-1}\circ\lim$. The following uniform version of the
low basis theorem of Jockusch and Soare~\cite{JS72} was proved by
the author, de Brecht and Pauly~\cite{BBP12}.

\begin{theorem}[Uniform low basis theorem]
\label{thm:ULB}
\label{thm:low-basis}
$\C_{2^\IN}\lW\C_\IR\lW\L\lW\lim$.
\end{theorem}

This implies, in particular, $\C_\IN\lW\L$.

There are many further studies, e.g., on Nash equilibria by Arno Pauly~\cite{Pau10,Pau11} that can also be described with the help of choice,
on probabilistic versions of choice, for instance, by the author and Pauly~\cite{BP10}, the author, Gherardi and H\"olzl~\cite{BGH15a},
Bienvenu and Kuyper~\cite{BK17}, on degrees defined by jumps of choice, e.g., by the author, Gherardi and Marcone~\cite{BGM12},
non-standard degrees obtained by Ramsey's theorem, by 
Dorais, Dzhafarov, Hirst, Mileti and Shafer~\cite{DDH+16},
Hirschfeldt and Jockusch~\cite{Hir15,HJ16},
the author and Rakotoniaina~\cite{BR17},
Dzhafarov, Goh, Hirschfeldt, Patey and Pauly~\cite{DGH+20},
Cholak, Dzhafarov, Hirschfeldt and Patey~\cite{CDHP20}, 
Marcone and Valenti~\cite{MV21}, and Dzhafarov and Patey~\cite{DP21}.

The classifications presented here are essentially in line with results from reverse mathematics~\cite{Sim99}.
There are several attempts to establish formal bridges between the proof-theoretic side of reverse mathematics
to the more uniform computational side of Weihrauch complexity.
These can be found, for instance, in work of Fujiwara~\cite{Fuj21}, Uftring~\cite{Uft21}, and Kuyper~\cite{Kuy17}.

\section{Degrees of Translators Between Real Number Representations}

We are now going to discuss some of Weihrauch's early results from \cite{Wei92a,Wei92c}.
We will try to put them into the context of more recent results with the aim to show that many important
Weihrauch degrees were already anticipated in these initial reports.
The main objective of \cite{Wei92a} is the study of a kind of quotient structure of the
diagram displayed in Figure~\ref{fig:real-number-representations}. For any two representations $\delta_1,\delta_2$
of some set $X$ we can consider the {\em implication problem} $(\delta_1\to\delta_2):\In\IN^\IN\mto\IN^\IN$
\[(\delta_1\to\delta_2)(p):=\{q\in\IN^\IN:\delta_2(q)=\delta_1(p)\}\]
with $\dom(\delta_1\to\delta_2):=\dom(\delta_1)$, which measures the complexity of translating $\delta_1$ into $\delta_2$.
In particular, $(\delta_1\to\delta_2)$ is computable (in the sense that it has a computable refinement) if and only
if $\delta_1\leq\delta_2$. A benchmark that can be used to measure the complexity of other problems
is the problem 
\[\EC:\IN^\IN\to2^\IN,p\mapsto\range(p-1)\] 
that translates an enumeration of a set into its characteristic 
function. 
Using this terminology one of the results obtained by Weihrauch is the following.

\begin{theorem}[Weihrauch~{\cite[Theorem~4]{Wei92a}}]
\label{thm:Weihrauch4}
We obtain
\[\EC\equivW(\rho'\to\rho)\equivW(\rho'\to\rho_>)\equivW(\rho_<\to\rho_>)\equivW(\rho_<\to\rho).\]
\end{theorem}

The proof of this proposition builds on earlier results by von Stein~\cite{Ste89}.
While Theorem~\ref{thm:Weihrauch4} describes quotients in the upper part of the diagram
in Figure~\ref{fig:real-number-representations}, the next result describes quotients that refer to the lower part of the diagram.
In order to capture the complexity, one needs to modify the benchmark problem $\EC$ as follows.
By $\EC_1:\In\IN^\IN\to2^\IN$ we denote the restriction of $\EC$ to such enumerations $p$ for which 
\[|\IN\setminus\range(p-1)|\leq1\]
holds, i.e., the enumerated set $\range(p-1)$ 
is either $\IN$ or $\IN\setminus\{k\}$ for some $k\in\IN$.

\begin{theorem}[Weihrauch~{\cite[Theorem~22]{Wei92a}}]
\label{thm:Weihrauch22}
We obtain
\[\EC_1\equivW(\rho\to\rho_\mathrm{cf})\equivW(\rho\to\rho_\mathrm{cf>})\equivW(\rho_\mathrm{cf<}\to\rho_\mathrm{cf>})\equivW(\rho_\mathrm{cf<}\to\rho_\mathrm{cf}).\]
\end{theorem}

The proof of this result can be built on the observation that $\EC_1$ is equivalent to the rationality problem
(see Proposition~\ref{prop:rationality} below), and that relative to the rationality problem the continued fraction algorithm
is computable. We mention one last result from \cite{Wei92a} that is related to the separation problem.
Similarly as with $\EC$, the separation problem $\SEP$
has a version $\SEP_1$ restricted to the separation of sets which are almost complementary in the sense that
the pairs $\langle p,q\rangle\in\dom(\SEP_1)$ have to satisfy the additional condition
\[|\IN\setminus(\range(p-1)\cup\range(q-1))|\leq1.\]
Using this terminology, we can say something on the middle part of the diagram in Figure~\ref{fig:real-number-representations}.
Namely, the problem $\SEP_1$ characterizes exactly the complexity of translating the Cauchy representation
into the decimal representation.

\begin{theorem}[Weihrauch~{\cite[Theorem~13]{Wei92a}}]
$\SEP_1\equivW(\rho\to\rho_{10}).$
\end{theorem}

There are further results along these lines included in \cite{Wei92a} that we cannot all summarize here.
The three Weihrauch degrees of $\EC$, $\EC_1$ and $\SEP_1$ have anti\-cipated
classes that have been studied much later. We mention some of these observations.

\begin{theorem}
\label{thm:EC-equivalents}
We obtain:
\begin{enumerate}
\item  $\SEP\equivW\WKL\leqW\EC\equivW\lim$, 
\item $\SEP_1\equivW\C_{\#\leq2}\leqW\EC_1\equivW\SORT$.
\end{enumerate}
\end{theorem}

The equivalence $\EC\equivW\lim$ follows from \cite[Proposition]{Bra05} (see also \cite[Lemma~5.3]{BG11} and \cite{BBP12}).
The equivalence $\EC_1\equivW\SORT$ is proved in Proposition~\ref{prop:rationality} in the appendix,
where also a definition of the problem $\SORT$ is given, which was originally introduced by Pauly and Neumann~\cite{NP18}.
As mentioned before, the equivalence $\SEP\equivW\WKL$ has been proved by Gherardi and Marcone \cite[Theorem~6.7]{GM09}.
The reduction $\WKL\leqW\lim$ follows 
from the uniform version of the low basis theorem (Theorem~\ref{thm:ULB}) or can easily be proved directly. 

Finally, $\C_{\#\leq2}$ is the closed choice problem on Cantor space restricted to sets with one or two elements.
This problem was introduced and studied extensively by Le Roux and Pauly~\cite{LRP15a}.
The proof of $\SEP_1\equivW\C_{\#\leq2}$ is sketched in Proposition~\ref{prop:sep-choice} in the appendix.
The fact $\SEP_1\leqW\EC_1$ is stated by Weihrauch  in \cite[Section~6]{Wei92a}.
In the same section we also find the following statement (the first separation as conjecture).

\begin{proposition}[Weihrauch {\cite[Section~6]{Wei92a}}]
$\SEP_1\lW\EC_1\lW\EC$.
\end{proposition}

Weihrauch attributes the separation in the second claim to Matthias Schr\"oder.
Indeed this follows from the fact that $\SORT$
always has a computable output, whereas $\EC$ can have non-computable outputs on computable inputs.
The other separation $\EC_1\nleqW\SEP_1$, mentioned as
conjecture by Weihrauch, can be proved as follows:
we have $\LPO\leqW\EC_1$, but $\LPO\nleqW\SEP_1$, since
otherwise $\lim\equivW\widehat{\LPO}\leqW\WKL$ would follow,
which is wrong by the low basis theorem (Theorem~\ref{thm:ULB}).
The details and notations used in this proof will be discussed in the next section.

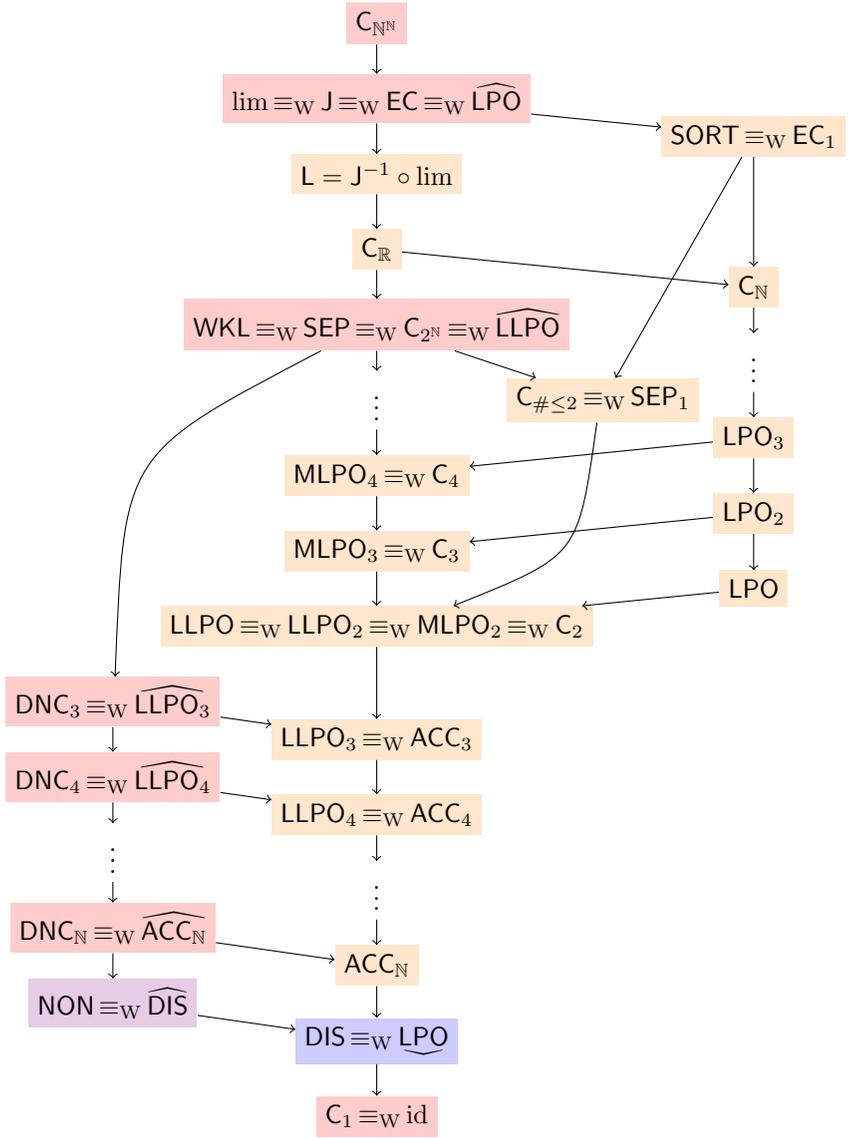
\begin{figure}[tb]
\begin{tikzpicture}[scale=1,auto=left,every node/.style={fill=gray!20}]
\node[style={fill=red!20}]  (C1) at (0,-1.5) {$\C_1\equivW\id$};
\node[style={fill=blue!20}]  (DIS) at (0,-0.5) {$\DIS\equivW\uwidehat{\LPO}$};
\node[style={fill=orange!20}]  (LLPO8) at (0,0.5) {$\ACC_\IN$};
\node[style={fill=gray!0}] (LLPOvdots) at (0,1.5) {$\vdots$};
\node[style={fill=orange!20}]  (LLPO4) at (0,2.5) {$\LLPO_4\equivW\ACC_4$};
\node[style={fill=orange!20}]  (LLPO3) at (0,3.5) {$\LLPO_3\equivW\ACC_3$};
\node[style={fill=orange!20}]  (LLPO) at (0,5) {$\LLPO\equivW\LLPO_2\equivW\MLPO_2\equivW\C_2$};
\node[style={fill=orange!20}]  (MLPO3) at (0,6) {$\MLPO_3\equivW\C_3$};
\node[style={fill=orange!20}]  (MLPO4) at (0,7) {$\MLPO_4\equivW\C_4$};
\node[style={fill=gray!0}] (MLPOvdots) at (0,8) {$\vdots$};
\node[style={fill=red!20}] (WKL) at (0,9) {$\WKL\equivW\SEP\equivW\C_{2^\IN}\equivW\widehat{\LLPO}$};
\node[style={fill=orange!20}] (L) at (0,11) {$\L=\J^{-1}\circ\lim$};
\node[style={fill=red!20}] (LIM) at (0,12) {$\lim\equivW\J\equivW\EC\equivW\widehat{\LPO}$};
\node[style={fill=orange!20}] (CR) at (0,10) {$\C_\IR$};
\node[style={fill=violet!20}]  (NON) at (-3.5,0) {$\NON\equivW\widehat{\DIS}$};
\node[style={fill=red!20}]  (DNCN) at (-3.5,1) {$\DNC_\IN\equivW\widehat{\ACC_\IN}$};
\node[style={fill=gray!0}] (DNCvdots) at (-3.5,2) {$\vdots$};
\node[style={fill=red!20}]  (DNC4) at (-3.5,3) {$\DNC_4\equivW\widehat{\LLPO_4}$};
\node[style={fill=red!20}]  (DNC3) at (-3.5,4) {$\DNC_3\equivW\widehat{\LLPO_3}$};
\node[style={fill=red!20}]  (CB) at (0,13) {$\C_{\IN^\IN}$};
\node[style={fill=orange!20}]  (LPO) at (5,5.5) {$\LPO$};
\node[style={fill=orange!20}]  (LPO3) at (5,6.5) {$\LPO_2$};
\node[style={fill=orange!20}]  (LPO4) at (5,7.5) {$\LPO_3$};
\node[style={fill=gray!0}] (LPOvdots) at (5,8.5) {$\vdots$};
\node[style={fill=orange!20}]  (SORT) at (5,11.5) {$\SORT\equivW\EC_1$};
\node[style={fill=orange!20}]  (CN) at (5,9.5) {$\C_\IN$};
\node[style={fill=orange!20}]  (SEP1) at (3,8) {$\C_{\#\leq2}\equivW\SEP_1$};

\draw[->] (CB) edge (LIM);
\draw[->] (LIM) edge (L);
\draw[->] (L) edge (CR);
\draw[->] (CR) edge (WKL);
\draw[->] (CR) edge (CN);
\draw[->] (WKL) edge (MLPOvdots);
\draw[->] (MLPOvdots) edge (MLPO4);
\draw[->] (MLPO4) edge (MLPO3);
\draw[->] (MLPO3) edge (LLPO);
\draw[->] (LLPO) edge (LLPO3);
\draw[->] (LLPO3) edge (LLPO4);
\draw[->] (LLPO4) edge (LLPOvdots);
\draw[->] (LLPOvdots) edge (LLPO8);
\draw[->] (LLPO8) edge (DIS);

\draw[->,bend right=30, looseness=1.5] (WKL) edge (DNC3);
\draw[->] (DNC3) edge (DNC4);
\draw[->] (DNC4) edge (DNCvdots);
\draw[->] (DNCvdots) edge (DNCN);
\draw[->] (DNC3) edge (LLPO3);
\draw[->] (DNC4) edge (LLPO4);
\draw[->] (DNCN) edge (LLPO8);
\draw[->] (DNCN) edge (NON);
\draw[->] (NON) edge (DIS);
\draw[->] (DIS) edge (C1);

\draw[->] (LIM) edge (SORT);
\draw[->] (SORT) edge (CN);
\draw[->] (CN) edge (LPOvdots);
\draw[->] (LPOvdots) edge (LPO4);
\draw[->] (LPO4) edge (LPO3);
\draw[->] (LPO3) edge (LPO);
\draw[->] (LPO) edge (LLPO);
\draw[->] (LPO3) edge (MLPO3);
\draw[->] (LPO4) edge (MLPO4);

\draw[->] (WKL) edge (SEP1);
\draw[->] (SORT) edge (SEP1);
\draw[->,bend left=30, looseness=1.5] (SEP1) edge (LLPO);
\end{tikzpicture}
\caption{Some Weihrauch degrees.}
\label{fig:Weihrauch-lattice}
\end{figure}

The diagram in Figure~\ref{fig:Weihrauch-lattice} shows a portion of the Weihrauch lattice
that includes the problems discussed here together with problems mentioned in the next section.
Any arrow indicates a Weihrauch reduction against the direction of the arrow (which is natural, as arrows correspond to logical implications in this way).
The diagram should be complete up to transitivity and reflexivity.

\clearpage

\section{Degrees of Omniscience Principles}

In constructive mathematics principles of omniscience represent 
certain unsolvable problems of different degrees of complexity.
Bishop introduced the {\em limited principle of omniscience} $\LPO$ and
the {\em lesser limited principle of omniscience} $\LLPO$ \cite{Bis67}.
Richman~\cite{Ric90} generalized the latter mentioned principle to
principles $\LLPO_n$. Weihrauch~\cite{Wei92c} translated these
principles into problems in his lattice structure. 
We formulate them as multi-valued maps for all $n\geq1$,
where $\langle \,\rangle$ denotes some standard tupling function on Baire space:

\begin{itemize}
\item $\LPO_n:\IN^\IN\to\IN,\langle p_1,...,p_n\rangle\mapsto|\{i\in\IN:p_i=\widehat{0}\}|$,
\item $\LLPO_n:\In\IN^\IN\mto\IN,\langle p_1,...,p_n\rangle\mapsto\{i\in\IN:p_i=\widehat{0}\}$ with\\
 $\dom(\LLPO_n):=\{\langle p_1,...,p_n\rangle:|\{i\in\IN:p_i\not=\widehat{0}\}|\leq1\}$,
\item $\MLPO_n:\In\IN^\IN\mto\IN,\langle p_1,...,p_n\rangle\mapsto\{i\in\IN:p_i=\widehat{0}\}$ with\\
 $\dom(\MLPO_n):=\{\langle p_1,...,p_n\rangle:|\{i\in\IN:p_i=\widehat{0}\}|\geq1\}$.
\end{itemize}

Using these general problems we define $\LPO:=\LPO_1$
and $\LLPO:=\LLPO_2=\MLPO_2$.
Among other results Weihrauch proved the following properties of these problems.

\begin{theorem}[Weihrauch~{\cite[Theorems~4.2, 4.3, 5.4]{Wei92c}}]
\label{thm:Weihrauch-omniscience}
For all $n,k\geq1$ we obtain
\begin{enumerate}
\item $\LPO_n\lW\LPO_{n+1}$,
\item $\LLPO_{n+2}\lW\LLPO_{n+1}\lW\LPO$,
\item $\MLPO_{n+1}\lW\MLPO_{n+2}$,
\item $\MLPO_{n+1}\lW\LPO_n$,
\item $\LPO_k\nleqW\MLPO_{n+1}$.
\end{enumerate}
\end{theorem}

The diagram in Figure~\ref{fig:Weihrauch-lattice} illustrates these results.
Later on, it was noticed that some of these omniscience proplems can
also be seen as choice problems. By $\C_n$ we denote closed choice
as defined in the previous section for the space $n=\{0,1,...,n-1\}$.
By $\ACC_n$ we denote the {\em all-or-co-unique choice problem}, which
is $\C_n$ restricted to sets $A\In n$ with $|n\setminus A|\leq1$.
This definition is used analogously for $\C_\IN$ and $\ACC_\IN$ 
and the space $\IN=\{0,1,2,...\}$ of natural numbers.
The following was observed in \cite[Example~3.2]{BBP12}
and \cite[Fact~3.3]{BHK17a}.

\begin{proposition}[Choice and omniscience]
\label{prop:choice-omniscience}
We obtain
$\MLPO_n\equivW\C_n$ and $\LLPO_n\equivW\ACC_n$ for all $n\geq2$.
\end{proposition}

Some of the reductions and separations shown in Figure~\ref{fig:Weihrauch-lattice}
have been proved elsewhere. We mention an example on choice and cardinality
proved by Le Roux and Pauly~\cite[Figure~1]{LRP15a}.

\begin{proposition}[Le Roux and Pauly~{\cite{LRP15a}}]
We obtain
\begin{enumerate}
\item $\C_2\lW\C_{\#\leq2}\nleqW\C_\IN$,
\item $\C_3\nleqW\C_{\#\leq2}$.
\end{enumerate}
\end{proposition}

Other interesting results are related to the parallelization of omniscience problems.
For any problem $f:\In X\mto Y$ 
\[\widehat{f}:\In X^\IN\mto Y^\IN,(x_n)_n\mapsto\bigtimes_{n\in\IN}f(x_n)\]
is called the {\em parallelization} of $f$.
This concept was introduced in \cite{BG11}, and it was proved that $f\mapsto\widehat{f}$
is a closure operator in the Weihrauch lattice. The following result characterizes the
parallelizations of $\LPO$ and $\LLPO$.

\begin{proposition}[B.\ and Gherardi~{\cite[Theorem~8.2, Lemma~6.3]{BG11}}]
We have:
$\widehat{\LLPO}\equivW\WKL$ and $\widehat{\LPO}\equivW\EC$.
\end{proposition}

It was discovered by Higuchi and Kihara~\cite{HK14a} 
that the parallelization of $\LLPO_n$
is related to the problems of {\em diagonally non-computable functions}
\[\DNC_X:\IN^\IN\mto X^\IN,p\mapsto\{f\in X^\IN:\mbox{$f$ is diagonally non-comp.\ relative to $p$}\},\]
defined for all $X\In\IN$.
As usual, $f:\IN\to X$ is called {\em diagonally non-computable relative to} $p$,
if $(\forall i)\;\varphi_i^p(i)\not=f(i)$, where $\varphi_i^p$ is a G\"odel numbering
of the functions that are computable relative to $p\in\IN^\IN$. We now obtain the following
parallelizations.

\begin{proposition}[Higuchi and Kihara {\cite[Proposition~81]{HK14a}} and B., Hendtlass and Kreuzer~{\cite[Theorem~5.2]{BHK17a}}]
$\DNC_n\equivW\widehat{\LLPO_n}$ for all $n\geq2$ and $\DNC_\IN\equivW\widehat{\ACC_\IN}$.
\end{proposition}

It is interesting to note that Jockusch~\cite{Joc89} proved in 1989 a statement
that is very similar to Weihrauch's result $\LLPO_{n+3}\lW\LLPO_{n+2}$ from Theorem~\ref{thm:Weihrauch-omniscience}
for the unparallelized case. Jockusch result was expressed in terms
of Medvedev reducibility, but the separations imply separations of the
corresponding Weihrauch degrees.

\begin{proposition}[Jockusch~{\cite[Theorem~6]{Joc89}}]
We obtain for all $n\in\IN$:
$\DNC_\IN\lW\DNC_{n+3}\lW\DNC_{n+2}$.
\end{proposition}

In fact, the results of Jockusch and Weihrauch can be jointly generalized to $\ACC_{n+2}\nleqW\DNC_{n+3}$~\cite[Proposition~5.7]{BHK17a}.

For completeness, we have also added to Figure~\ref{fig:Weihrauch-lattice} the lowest known natural discontinuous Weihrauch
degree, namely that of the {\em discontinuity problem} $\DIS$~\cite{Bra22}.
This problem can be defined by {\em stashing} $\uwidehat{\LPO}$ of $\LPO$~\cite{Bra21}  
(stashing is an operation dual to parallelization).
We mention that the {\em non-computability problem} $\NON$, which
was introduced and studied in \cite{BHK17a}, and maps every $p\in\IN^\IN$
to some $q\in\IN^\IN$ that is not computable from $p$, is the
parallelization of $\DIS$~\cite{Bra21}.
We omit the definitions here and point the interested reader to the references.

The main purpose of the diagram in Figure~\ref{fig:Weihrauch-lattice}
is to illustrate that some of the problems discussed by Weihrauch in his
initial two reports are related to important cornerstones among the Weihrauch degrees, some of which were much later rediscovered and studied under different names.  All parallelizable and stashable degrees in Figure~\ref{fig:Weihrauch-lattice} are displayed in violet, all only parallelizable degrees are displayed in red, all only stashable degrees in blue, and all other degrees in orange.

\section{Some Work by Students in Hagen}

The purpose of this section is to briefly discuss some follow-up work on Weihrauch's initial reports~\cite{Wei92a,Wei92c}.
In fact, Klaus Weihrauch supervised six MSc and PhD theses 
on topics related to his reducibility over a period of 18 years.
The first of these even appeared before \cite{Wei92a,Wei92c} and are listed here
in chronological order. 
Like \cite{Wei92a,Wei92c} most of this work remained unpublished for quite some
time, and we give some
references to later publications that overlap with the content of these theses.

\begin{itemize}
\item {\em Diploma thesis of Torsten von Stein}~\cite{Ste89}:
in this thesis the concept of Weihrauch reducibility appears in writing for the first time
and the author introduces the so-called C-hierarchy that measures the complexity
of a problem $f$ by the number of applications of the problem $\EC$ that are required
to solve $f$. He also studies quite a number of specific
problems from analysis, most notably the differentiation problem (see Theorem~\ref{thm:von-Stein} below).
The results remain unpublished except for some pointers to this work elsewhere. 
\item {\em Diploma thesis of Uwe Mylatz}~\cite{Myl92}:
Uwe Mylatz continued the work of von Stein, and he studied also the second and third level of the C-hierarchy.
Among other problems, he looked at higher derivatives, rationality of reals, monotonicity of sequences,  convergence, density etc.
This work remained unpublished.
\item {\em Diploma thesis of the author of these notes}~\cite{Bra93}:
in this thesis it is essentially proved that the C-hierarchy is identical to the effective
Borel hierarchy. 
A much extended version of this thesis was later published in \cite{Bra05}, which also
includes a proof that for admissible representations of computable Polish spaces
Borel measurability on a certain finite level corresponds to the corresponding measurability of the realizers.
\item {\em PhD thesis of Peter Hertling}~\cite{Her96d}:
the author introduces methods of combinatorial nature that allow to classify the (strong) topological Weihrauch
degrees of certain problems with finite or countable discrete image, using forests of trees. These trees essentially capture
the nature of the occurring discontinuities. Some of these results were already touched upon
in the technical reports \cite{Her93,Her93a,Her93b,HW94a} and the proceedings article \cite{HW94}.
A significant extension of some of the crucial results in \cite{Her96d} was finally published in \cite{Her20}.
This thesis has also influenced later work on this topic, for instance \cite{Sel07,Sel11,HS14a,KM19}.
\item {\em PhD thesis of Uwe Mylatz}~\cite{Myl06}:
in his PhD thesis the author continues the study of the problems $\LPO_n, \LLPO_n, \MLPO_n$ and many variants thereof
of combinatorial nature. Partially the methods developed in \cite{Her96d} are considered for separation results.
The thesis remains unpublished.
\item {\em MSc thesis of Arno Pauly}~\cite{Pau07}: in this thesis the author studies the classification of a number
of concrete problems and the structure of the lattice induced by the continuous version of Weihrauch reducibility from a lattice-theoretic perspective. 
The publication~\cite{Pau10a} extends many of the structural results.
\end{itemize}

Presumably, the first time Weihrauch reducibility appeared in a published article was in~\cite{Bra99},
where the author generalized Pour-El and Richards first main theorem~\cite{PR89}.
For instance the fact that the operator $d$ of differentiation is linear, has a closed graph and satisfies 
some minimal computability properties implies already the direction $\EC\leqW d$ in the following result.

\begin{theorem}[von Stein~{\cite[Theorem~12]{Ste89}}]
\label{thm:von-Stein}
$\EC\equivW d$ where $d$ is the operator of differentiation $d:\In\CC[0,1]\to\CC[0,1],f\mapsto f'$, restricted to
continuously differentiable functions.
\end{theorem}

\section{Conclusion}

The purpose of these notes is to show how the subject of Weihrauch complexity emerged in work
of the Hagen school of computable analysis, in particular, from Weihrauch's original contributions
to this topic and work produced by his students. Many more recent developments and prominent
degrees were already anticipated in this early work.

Today Weihrauch complexity is a mature research topic that has found interest from researchers
in computability theory, computable analysis, reverse mathematics, and proof theory.
This is witnessed by two Dagstuhl seminars in 2015 and 2018 that were dedicated to Weihrauch complexity
and a number of international PhD and MSc theses that are either mostly
or partially dedicated to studying Weihrauch complexity.
This includes the theses by  Gherardi~\cite{Ghe06b}, Pauly~\cite{Pau11}, Higuchi~\cite{Hig12a}, Carroy~\cite{Car13},
Neumann~\cite{Neu14,Neu18}, 
Rakotoniaina~\cite{Rak15}, Borges~\cite{Bor16}, Patey~\cite{Pat16}, Sovine~\cite{Sov17}, Nobrega~\cite{Nob18},
 Thies~\cite{Thi18}, Uftring~\cite{Uft18}, Goh~\cite{Goh19a}, Angl{\'e}s d'Auriac~\cite{AdA20}, and Valenti~\cite{Val21} 
 (in chronological order from 2011 onwards, not mentioning the theses discussed in the previous section).
A comprehensive up-to-date bibliography on Weihrauch complexity can be found online.~\footnote{\tt http://cca-net.de/publications/weibib.php}

\section*{Appendix: Some Proofs}

In this appendix we add some proofs that justify some claims made earlier.
It seems that most of the corresponding equivalences stated here have not been noticed
widely and do not appear in writing elsewhere.

We recall that the {\em sorting problem} $\SORT:\{0,1\}^\IN\to\{0,1\}^\IN$
is defined by 
\[\SORT(p):=\left\{\begin{array}{ll}
0^n\widehat{1} & \mbox{if $p$ contains exactly $n$ zeros}\\
\widehat{0}      & \mbox{if $p$ contains infinitely many zeros}
\end{array}\right.\]
This problem was introduced by Pauly and Neumann~\cite{NP18}.
Here $\widehat{n}=nnn...\in\IN^\IN$ denotes the constant sequence with value $n\in\IN$.

We define the {\em rationality problem} $\RAT:\IR\mto\{0,1\}^\IN$ by
\[\RAT(x):=\left\{\begin{array}{ll}
\{0^n\widehat{1}:x=q_n\} & \mbox{if $x\in\IQ$}\\
\{\widehat{0}\} & \mbox{if $x\in\IR\setminus\IQ$}
\end{array}\right.,\]
where $(q_n)_{n\in\IN}$ is the standard enumeration of the 
rational numbers $\IQ$ given by $q_{\langle n,k,m\rangle}:=\frac{n-k}{m+1}$, where $\langle\,\rangle$ denotes
a standard Cantor tupling on $\IN$.
That is, $\RAT$ makes the property of being a rational number c.e.\ and on top
of it, it determines the value of the rational number in the case of a rational input $x\in\IR$.
For irrational inputs $x$ no information is provided.
Since the representations of the reals below and including $\rho$ in the diagram in Figure~\ref{fig:real-number-representations}
are all computably equivalent, if restricted to irrational numbers, it is exactly the problem $\RAT$ that
can help to translate the representations into each other. We prove the following statement about
the equivalence class of $\RAT$.

\begin{proposition}
\label{prop:rationality}
$\EC_1\equivW\RAT\equivW\SORT$.
\end{proposition}
\begin{proof}
We are going to prove $\RAT\leqSW\EC_1\leqSW\SORT\leqW\RAT$.

We start with proving $\RAT\leqSW\EC_1$. 
Given an input $x\in\IR$ for $\RAT$, we start producing an enumeration $p\in\IN^\IN$ 
of $\IN\setminus\{0\}=\range(p-1)$ on the output side.
Simultaneously, we try falsifying $x=q_0$. 
If it turns out that $x\not=q_0$, then we enumerate $0$ into our set and we 
switch to producing an enumeration of $\IN\setminus\{m\}$
with some sufficiently large $m\in\IN$ that has not yet been enumerated and such that $q_m=q_1$.
Simultaneously, we try falsifying $x=q_m$. 
If we continue inductively in this way with some $m$ with $q_m=q_2$ in the next step, 
then for rational input $x\in\IQ$ we produce some enumeration $p$
of a set $A=\IN\setminus\{m\}=\range(p-1)$ with $q_m=x$ and for irrational $x$ we produce some enumeration $p$
of $A=\IN=\range(p-1)$.
Hence, given the characteristic function $\chi_A=\EC_1(p)$ we can compute a point in $\RAT(x)$ as follows:
if $\chi_A=\widehat{1}$, then $\widehat{0}\in\RAT(x)$ and if $\chi_A=1^{m}0\widehat{1}$, then 
$0^m\widehat{1}\in\RAT(x)$.

We now prove $\EC_1\leqSW\SORT$.
Given $p\in\IN^\IN$ that enumerates a set $A=\range(p-1)$ with $A=\IN$ or $A=\IN\setminus\{m\}$ for some $m\in\IN$,
we compute a sequence $q\in\{0,1\}^\IN$ as follows. Whenever we have found a consecutive
segment $\{0,1,....,k\}\In A$ of length $k+1$ in the enumeration $p$, then
we ensure that the output $q$ contains exactly $k+1$ zeros. As long as no other
information is available on the input side, we append digits $1$ to the output $q$. 
This algorithm ensures that $q$ contains infinitely many zeros (and hence $\SORT(q)=\widehat{0}$)
if and only if $A=\IN$
and $\SORT(q)=0^{m}111...$ if and only if $A=\IN\setminus\{m\}$.
Hence given $\SORT(q)$ we can compute the characteristic function $\chi_A$ of $A$ as follows:
if $\SORT(q)=\widehat{0}$ then $\chi_A=\widehat{1}$ and if $\SORT(q)=0^m\widehat{1}$ then 
$\chi_A=1^m0\widehat{1}$.

Now we prove $\SORT\leqW\RAT$. Given an input $p\in\{0,1\}^\IN$ we start producing
better and better approximations of $x_0=\frac{1}{2^{\langle 0,0\rangle}}\in\IR$ on the output side.
As soon as we find the first digit $0$ in $p$, we switch to writing better and better
approximations of some number $x_1=\frac{2m+1}{2^{\langle 1,k\rangle}}\in\IR$ on the output side with suitable $m,k\in\IN$,
such that $x_1$ is compatible with the previous approximations of $x_0$ and we start with producing an approximation of $x_1$
that is good enough such that no rational number with denominator smaller than 
$2^{\langle 1,k\rangle}$ can
satisfy it (which is possible as $2m+1$ and $2^{\langle 1,k\rangle}$ are coprime).
In general, if we find the $(n+1)$--th digit $0$ in $p$, then we 
switch to producing approximations of some number $x_{n+1}=\frac{2m+1}{2^{\langle n+1,k\rangle}}$ for appropriate $m,k\in\IN$, such that
$x_{n+1}$ is compatible with the previous approximation of $x_n$.
Again the first approximation of $x_{n+1}$ that we produce is good enough such that
no rational number with denominator smaller than $2^{\langle n+1,k\rangle}$ can satisfy it. 
If there are infinitely many zeros in $p$, then the output converges to $x:=\lim_{n\to\infty}x_n$,
which must be irrational, since larger and larger denominators are excluded by the conditions above.
If there is only a finite number $j$ of zeros in $p$, then the final output is $x=x_{j}$ and from $0^i\widehat{1}\in\RAT(x)$
we can compute $m,n,k\in\IN$ such that $x=q_i=\frac{2m+1}{2^{\langle n,k\rangle}}$, where the numbers
$m,n,k$ are unique since the numerator and denominator are coprime. Given $x$ and the original input $p$,
we can now compute $\SORT(p)$ as follows: we write as many zeros to the output as we can find in $p$.
Simultaneously, we try to find some $i$ with $0^i\widehat{1}\in\RAT(x)$. If we find such an $i$, then
we determine $n$ as above and we extend the output to $0^n\widehat{1}=\SORT(p)$. 
If we never find such an $i$, then $p$ contains infinitely many zeros and the output is $\widehat{0}=\SORT(p)$.
\end{proof}

The equivalence $\RAT\equivW\SORT$ is similar to the statement of \cite[Theorem~23]{NP18}.\footnote{However, 
\cite[Theorem~23]{NP18} does not seem to hold in full generality as stated, and the proof given here fills a gap and is an alternative proof for the case $X=\IR$.}

\begin{proposition}[Separation and choice]
\label{prop:sep-choice}
$\C_{\#\leq2}\equivSW\SEP_1$.
\end{proposition}
\begin{proof} (Sketch)
Firstly, we note that 
\[\dom(\SEP_1)=\{\langle p,q\rangle\in\IN^\IN:1\leq|\SEP\langle p,q\rangle|\leq2\}\]
and negative information on the closed set $\SEP_1\langle p,q\rangle\In2^\IN$ can be computed from $p,q$.
This provides the reduction $\SEP_1\leqSW\C_{\#\leq2}$. 

For the inverse reduction, we have given by negative information a closed set $A\In2^\IN$ with one or two points.
We can compute from the negative information an infinite binary tree $T$ with $[T]=A$.
This tree has one or two infinite paths.
We now follow the construction in the proof of~\cite[Theorem~6.7]{GM09}.
There a predicate $\varphi_T(s,i)$ is used, which indicates that for $s\in\{0,1\}^*$ and $i\in\{0,1\}$
there is a finite length $n\in\IN$ such that the path $si$ can be extended in $T$ up to length $n$, 
but not the path $s(i-1)$. 
Without loss of generality, we assume that words $s\in\{0,1\}^*$ are identified with numbers $s\in\IN$ that encode them with respect to some standard encoding.
In the proof of~\cite[Theorem~6.7]{GM09} it is shown how to compute $p,q\in\IN^\IN$ such that
\begin{itemize}
\item $\range(p-1)=\{s+2:\varphi_T(s,0)\}\cup\{0\}$,
\item $\range(q-1)=\{s+2:\varphi_T(s,1)\}\cup\{1\}$.
\end{itemize}
Since our tree $T$ has exactly one or two infinite paths, it follows that
$\varphi_T(s,0)\vee\varphi_T(s,1)$ holds for all $s$, except possibly one.
That is 
\[|\IN\setminus(\range(p-1)\cup\range(q-1))|\leq1.\] 
Hence, $\langle p,q\rangle\in\dom(\SEP_1)$, and 
the proof of the reduction $\WKL\leqSW\SEP$ in~\cite[Theorem~6.7]{GM09} extends to a proof of the reduction 
$\C_{\#\leq2}\leqSW\SEP_1$.
\end{proof}

\bibliographystyle{abbrv}
\bibliography{C:/Users/\input{"|kpsewhich -var-value USERNAME"}/Documents/Spaces/Research/Bibliography/lit}

\section*{Acknowledgments}

The author thanks Peter Hertling, Arno Pauly and Klaus Weihrauch for helpful comments on this article.
This work has been supported by the {\em National Research Foundation of South Africa} (Grant Number 115269).

\end{document}